\documentclass{article}
\usepackage[centertags]{amsmath}
\usepackage{amsfonts}
\usepackage{amssymb}
\usepackage{amsthm}
\usepackage{latexsym}
\usepackage{amsmath}
\usepackage{amsfonts}
\usepackage[margin=1in]{geometry}
\usepackage{graphicx}

\title{Increasing stability of the continuation for general elliptic equations of second order.}

\author{Victor Isakov}

\newtheorem{theorem}{Theorem}[section]

\newtheorem{corollary}[theorem]{Corollary}

\begin{document}

\maketitle

\section{Introduction}

The Cauchy Problem (or equivalently the continuation of solutions) for partial differential equations has a long and rich history, starting with the Holmgren-John theorem on uniqueness for equations with analytic coefficients. It is of great importance in the theories of  boundary control and of inverse problems.  In 1938 T. Carleman introduced a special exponentially weighted energy (Carleman type) estimates to handle non analytic coefficients. These estimates imply in addition some conditional H\"older type  bounds for solutions of this problem.  In 1960  \cite{J} F. John  showed that for the continuation of solutions to the Helmholtz equation from inside of the unit disk onto any larger disk the stability estimate which is uniform with respect to the wave numbers is still of logarithmic type. Logarithmic stability is quite damaging for numerical solution of many inverse problems. In recent papers \cite{AI}, \cite{AI1}, \cite{HI}, \cite{I2}, \cite{I3}, \cite{IK} it was shown that in a certain sense stability is always improving for larger $k$ under (pseudo) convexity conditions on the geometry of the domain and of the coefficients of the elliptic equation. These assumptions are not satisfied in the John's example.

In this paper we attempt to eliminate any convexity type condition on the elliptic
operator or the domain. Due to  the John' counterexample, one can not expect increasing
stability for all solutions without suitable a priori constraint. We show that (near Lipschitz) stability holds on a subspace of ("low frequency") solutions which is is growing with the wave number $k$ under some mild boundedness constraints on complementary "high frequency" part. Using new logarithmic stability estimates in the Cauchy problem for hyperbolic equations without any convexity assumptions we bound the norm of the "high frequency" part norm  and obtain a conditional stability bound  which is improving for larger wave numbers $k$.

 We  consider the Cauchy problem
\begin{equation}
\label{PD}
(A-i a_0k+k^2) u= f\;\mbox{in}\;\Omega,
\end{equation}
with the Cauchy data
\begin{equation}
{\label{CauchyD}}
   u= u_{0},  \partial_{\nu} u=u_{1}\;  \text{ on }\; \Gamma\subset \partial \Omega,
\end{equation}
where
$$
Au=\sum_{j,m=1}^n a_{jm}\partial_j\partial_m u+
\sum_{j=1}^n a_{j}\partial_j u+ au
$$
is a general partial differential operator of second order satisfying the uniform ellipticity condition
$$
\varepsilon_0|\xi|^2\leq \sum_{j,l=1}^n a_{jl}(x)\xi_j\xi_l\leq E_0^2|\xi|^2
$$
for some positive numbers $\varepsilon_0, E_0$ and all $x\in\Omega$ and $\xi\in \mathbb{R}^n$.  
We assume that $a_{jm}, \partial_p a_{jm}, a_j, a\in L^{\infty}(\Omega)$.   

We consider bounded open domain $\Omega$ in ${\mathbb R}^n$ with the boundary
$\partial \Omega \in C^2$ and $\Gamma$ which is an open subset of $\partial\Omega$.

We use the classical Sobolev spaces $H^{(l)}(\Omega)$ with the standard norm
$\|\cdot\|_{(l)}(\Omega)$.
We recall the common notation
$$
\|u\|_{L^2((T,T);H^1(\Omega))}=
(\int_{-T}^{T} \|u(t , )\|^2_{(1)}(\Omega) dt)^{\frac{1}{2}}.
$$

In what follows $C$ denote generic constants depending only on $\Omega, \Omega_0,
\Gamma, \omega, A, a_0$.

 We   claim 

\begin{theorem}
Let $\Gamma$ be not empty set. Then there is a constant $C$ such that
\begin{equation}
\label{stability0}
\| u\|_{(0)}(\Omega)\leq
C(\|u\|_{(0)}(\Gamma)+k^{-1}(\|f\|_{(0)}(\Omega)+
\|u_0\|_{(1)}(\Gamma)+\|u_1\|_{(0)}(\Gamma)+\|u\|_{(1)}(\Omega)))
\end{equation}
 for all $ u\in H^{(2)}(\Omega)$ solving \eqref{PD}, \eqref{CauchyD}.
 
 Moreover, if in addition $\Omega$ is $C^2$ diffeomorphic to the unit ball, 
then there are a monotone family of closed subspaces $H^{(1)}(\Omega;k)$ of $H^{(1)}(\Omega)$ with $\overline{\cup_k H^{(1)}(\Omega;k)}=H^{(1)}(\Omega)$, 
 a semi norm 
 $|||  \cdot |||_{(1;k)}(\Omega)$ on $H^{(1)}(\Omega)$ which is zero on $H^{(1)}(\Omega;k)$ and decreasing with respect to $k$, and a constant $C$ such that
\begin{equation}
\label{stability1}
\| u\|_{(0)}(\Omega)\leq
C(\|u\|_{(0)}(\Gamma)+k^{-1}(\|f\|_{(0)}(\Omega)+
\|u_0\|_{(1)}(\Gamma)+\|u_1\|_{(0)}(\Gamma)+|||u|||_{(1,k)}(\Omega)))
\end{equation}
 for all $ u\in H^{(2)}(\Omega)$ solving \eqref{PD}, \eqref{CauchyD}. 
 \end{theorem}
 
 Now we give a version of this result for an increasing stability of the continuation from an open subset $\omega$ of $\Omega$.
 
 \begin{theorem}
Let $\omega$ be non empty open subset of $\Omega$.
Then there 
is a constant $C$ such that
\begin{equation}
\label{stability0om}
\| u\|_{(0)}(\Omega)\leq
C k^{-1}(\|f\|_{(0)}(\Omega)+
\|u\|_{(1)}(\omega)+||u||_{(1)}(\Omega))
\end{equation}
 for all $ u\in H^{(2)}(\Omega)$ solving \eqref{PD}.
 
 Moreover, if $\Omega$ is a 
 $C^2$-diffeomorphic to the unit ball and $\bar \omega$
 contains an open (in $\partial\Omega$) non empty subset of $\partial\Omega$, then there are a monotone family of closed subspaces $H^{(1)}(\Omega;k)$ of $H^{(1)}(\Omega)$ with $\overline{\cup_k H^{(1)}(\Omega;k)}=H^{(1)}(\Omega)$, 
 a semi norm 
 $|||  \cdot |||_{(1;k)}(\Omega)$ on $H^{(1)}(\Omega)$ which is zero on $H^{(1)}(\Omega;k)$ and decreasing with respect to $k$, and a constant $C$ such that
 \begin{equation}
\label{stability1om}
 \| u\|_{(0)}(\Omega)\leq
C(k^{-1}(\|f\|_{(0)}(\Omega)+
\|u\|_{(1)}(\omega)+|||u|||_{(1,k)}(\Omega)))
\end{equation}
 for all $ u\in H^{(2)}(\Omega)$ solving \eqref{PD}. 
 \end{theorem}

  Now, for particular $\Omega, \Gamma $ in ${\mathbb R}^n, n=2,3,$ we will describe more explicitly
 the subspaces $H^1(\Omega;k)$ of Lipschitz stability and semi norms $||| u |||_{(1,k)}$.   Let $\Omega=\{x:1<|x|<R\}$, $\Gamma=\{x:|x|=1\}$,
 and $\Gamma_1= \{x:|x|=R\}$.
  We will use the (angular) orthogonal series 
$$
u(,\sigma)=\sum_{1 \geq m, p\leq p(m)} u(m,p)e(\sigma;m,p),
$$
where $\{e(\sigma;m,p)\}, m=1,2,...,p=1,...,p(m)$ is an orthonormal basis of exponential functions or of spherical harmonics in $L^2(S^n)$. We recall that in the polar coordinates $(r,\sigma), \sigma \in S^n$ the Laplace operator
$\Delta=(\partial_r)^2+
r^{-2}\Delta_{\sigma}+
2r^{-1}\partial_r$, where 
$\Delta_{\sigma}$ is the Beltrami operator on the unit sphere $S^n$
in $\mathbb R^n$. As known,
$e(\sigma;m,p)$ are eigen functions of the Beltrami operator. 
If $n=2$, then corresponding 
 eigenvalues are
$\lambda_m=-(m-1)^2$ and $p(1)=1, p(m)=2, m=2,3,...$. if $n=3$ eigenvalues are $\lambda_m=-m(m-1)$ 
and $p(m)=2m-1$.
We  introduce the low frequency part of $u$
  \begin{equation}
   \label{ulp}
    u_l(,\sigma)=\sum_{E_0^2m^2<
    (1-\varepsilon)k^2, 1\leq p \leq p(m)} u(m,p)e(\sigma;m,p).
   \end{equation}

Under a constraint on the high frequency component of $u$ we have 
 
\begin{theorem}

Let $\theta>0$.

There are $C, C(\theta)$ such that for a solution
to the Cauchy problem
\eqref{PD}, \eqref{CauchyD}
$$
k\| u\|_{(0)}(\Gamma_1)+\|\nabla u\|_{(0)}(\Gamma_1)+\| u\|_{(1)}(\Omega)\leq 
$$
\begin{equation}
\label{stability3}
C(k\|u_0\|_{(0)}(\Gamma)+\|u_1\|_{(0)}(\Gamma)+\|f\|_{(0)}(\Omega))+
C(\theta)k^{-\frac{1}{2}+\theta}\|u-u_l\|_{(2)}(\Omega).
\end{equation}
Moreover,
\begin{equation}
\label{stability4}
\| u\|_{(1)}(\Omega)
\leq C(k\|u_0\|_{(0)}(\Gamma)+\|u_1\|_{(0)}(\Gamma)+\|f\|_{(0)}(\Omega)+
k^{-1}\|u-u_l\|_{(2)}(\Omega)),
\end{equation}
and
\begin{equation}
\label{stability5}
\| u\|_{(0)}(\Omega)\leq
C(\|u_0\|_{(0)}(\Gamma)+
k^{-1}(\|u_1\|_{(0)}(\Gamma)+\|f\|_{(0)}(\Omega)+
\|u-u_l\|_{(1)}(\Omega))).
\end{equation}
 \end{theorem}

Several ingredients of a proof can be found in \cite{I5}. Observe, that Theorem 1.1 implies a (best possible) Lipschitz stability when $u\in H^{1}(\Omega;k)$, and, 
since the subspaces $H^{1}(\Omega;k)$ grow with respect to the wave number $k$ and exhaust the whole $H^1$, can be viewed as an indication of the increasing stability 
in the Cauchy problem.

To bound the semi norm $|||u|||_{(1,k)}(\Omega)$  we will use a  logarithmic stability result for the hyperbolic equation
\begin{equation}
\label{PDH}
(A-a_0\partial_t -\partial_t^2) v= f_v\;\mbox{in}\; Q=\Omega\times(-T,T),
\end{equation}
with the Cauchy data
\begin{equation}
{\label{CauchyH}}
   v= v_0,  \partial_{\nu} v=v_1\;  \text{ on }\; \Gamma\times(-T,T),
\end{equation}
where $\Gamma$ is an open subset of $\partial\Omega$. Sometimes it is more convenient
to replace the Cauchy data
(\ref{CauchyH}) by $u$ on
$\omega\times(-T,T)$ where
$\omega=\Omega\cap V$, $V$ is a neighbourhood of a point of $\partial\Omega$. One of $\Gamma, \omega$ (but not the  both) can be the empty set.
 
\begin{theorem}

Let open $\Omega_0\subset \Omega,\bar{\Omega}_0\subset \Omega \cup \Gamma  $ and $0<T_0$.

Then there  are positive constants $C, \kappa, T$ depending only on $A,a_0,\Omega, \Omega_0, \Gamma, T_0$
  such that
  \begin{equation}
\label{stability2space}
\| \partial^{\alpha}v\|_{(0)}(\Omega_0\times (-T_0, T_0))\leq  C ((\|v\|_{(1)}(Q)+\|\partial_t \partial^{\alpha}v\|_{(0)}(Q))\frac{1}{|ln \varepsilon |}+
\varepsilon^{\kappa}\|v\|_{L^2((-T,T);H^1(\Omega))}),
\end{equation}
when $|\alpha|\leq 1, \alpha_0=0$ and
\begin{equation}
\label{stability2}
\| v\|_{(1)}(\Omega_0\times (-T_0, T_0))\leq  C ((\|v\|_{(1)}(Q)+\|\partial_t v\|_{(1)}(Q))\frac{1}{|ln \varepsilon |}+\varepsilon^{\kappa}
\|v\|_{L^2((-T,T);H^1(\Omega)})
\end{equation}
 for all $ v \in H^{(2)}(Q) $ solving \eqref{PDH}, \eqref{CauchyH}, 
where 
$$
\varepsilon =\|f_v\|_{(0)}(Q)+
\|v_0\|_{L^2((-T,T);H^1(\Gamma))} +
\|v_1\|_{(0)}(\Gamma\times(-T,T))+
\|v\|_{L^2((-T,T);H^1(\omega))}.$$ 
 \end{theorem}
 
 \begin{corollary}
 
 Let one of $\Gamma, \omega$ be not empty, $\Omega_0$ be a sub domain of $\Omega$ with $\bar\Omega_0\subset
 \Omega\cup \Gamma $.
 Assume that $\Omega$ is $C^2$-diffeomorphic to the unit ball.

There are a monotone family of closed subspaces $H^{(1)}(\Omega;k)$ of $H^{(1)}(\Omega)$ with $\overline{\cup_k H^{(1)}(\Omega;k)}=H^{(1)}(\Omega)$, a semi norm $|||  \cdot |||_{(1;k)}(\Omega)$ on $H^{(1)}(\Omega)$ which is zero on $H^{(1)}(\Omega;k)$ and decreasing with respect to $k$ 
and positive constants $C, \kappa$ such that 
$$
\| u\|_{(0)}(\Omega_0)\leq
C(\|u\|_{(0)}(\Gamma)+k^{-1}(\|f\|_{(0)}(\Omega)+
\|u_1\|_{(0)}(\Gamma)+\|u\|_{(1)}(\omega)+\varepsilon_1^{\kappa}\|u\|_{(1)}(\Omega))+
$$
\begin{equation}
\label{stability1E}
|||u|||^{\frac{1}{2}}_{(1,k)}(\Omega)
k^{-\frac{1}{2}}||u||^{\frac{1}{2}}_{(1)}(\Omega)\frac{1}{\sqrt{|ln \varepsilon_1|}} )
\end{equation}
 for all $ u\in H^{(2)}(\Omega)$ solving \eqref{PD}, \eqref{CauchyD}. Here
 $$
 \varepsilon_1=\|f\|_{(0)}(\Omega)+\|u_0\|_{(1)}(\Gamma)+
 \|u_1\|_{(0)}(\Gamma)+
 \|u\|_{(1)}(\omega).
 $$ 
 \end{corollary}

\begin{corollary}
 
 Let one of $\Gamma, \omega$ be not empty, $\Omega_0$ be a sub domain of $\Omega$ with $\bar\Omega_0\subset
 \Omega\cup \Gamma $.

There are  positive constants $C, \kappa$ such that 
$$
\| u\|_{(0)}(\Omega_0)\leq
C(\|u\|_{(0)}(\Gamma)+k^{-1}(\|f\|_{(0)}(\Omega)+
\|u_1\|_{(0)}(\Gamma)+\|u\|_{(1)}(\omega)+\varepsilon_1^{\kappa}\|u\|_{(1)}(\Omega))+
$$
\begin{equation}
\label{stability0E}
k^{-\frac{1}{2}}||u||_{(1)}(\Omega)\frac{1}{\sqrt{|ln \varepsilon_1|}} )
\end{equation}
 for all $ u\in H^{(2)}(\Omega)$ solving \eqref{PD}, \eqref{CauchyD}. Here
 $$
 \varepsilon_1=\|f\|_{(0)}(\Omega)+\|u_0\|_{(1)}(\Gamma)+
 \|u_1\|_{(0)}(\Gamma)+
 \|u\|_{(1)}(\omega).
 $$ 
 \end{corollary}

To derive Corollaries 1.5, 1.6 from Theorems 1.1, 1.2, and 1.4 we let $v(x,t)=u(x) e^{ikt}$. Then
$f_v(x,t)=f(x)e^{ikt}$,
$$
\|u\|_{(1)}(\Omega_0)\leq
C\sum_{|\alpha|\leq 1, \alpha_0=0}\| \partial^{\alpha}v\|_{(0)}(\Omega_0\times(-T_0,T_0)),
$$
$$
\|v\|_{(1)}(Q) +
\|\partial^{\alpha}\partial_t v\|_{(0)}(Q)\leq C k \|u\|_{(1)}(\Omega) \;, |\alpha|\leq 1, \alpha_0=0,
$$
$$
\|v_0\|_{L^2((-T,T);H^1(\Gamma))}\leq
C\| u_0\|_{(1)}(\Gamma),
\|v_1\|_{(0)}(\Gamma\times(-T,T))\leq
C\| u_1\|_{(0)}(\Gamma),
$$
$$
\|v\|_{L^2((-T,T);H^1(\omega))}\leq
C\| u\|_{(1)}(\omega),
\|v\|_{L^2((-T,T);H^1(\Omega))}\leq
C\| u\|_{(1)}(\Omega)
$$
and \eqref{stability1E},
\eqref{stability0E} follow
from \eqref{stability2space} and \eqref{stability0},\eqref{stability1}, \eqref{stability0om},\eqref{stability1om} used for $\Omega=\Omega_0$.

To illustrate that these results (in particular,
\eqref{stability1E},
\eqref{stability0E})  are sharp we recall the famous counterexample of Fritz John \cite{J}, \cite{I1}. Let $(r,\phi)$ be the polar coordinates in
$\mathbb R^2$. It is shown in \cite{J} that the functions $u(r,\phi;k)=
k^{-\frac{2}{3}}J_k(kr)e^{ik\phi}$ solve the homogeneous Helmholtz equation in the plane,
\begin{equation}
\label{J1}
k^{-1} \leq C\|u( ;k)\|_{(0)},\;
\|u( ;k)\|_{(1)}(\Omega)\leq C, \;
\|u( ;k)\|_{(1)}(\Gamma)+
\|\partial_r u( ;k)\|_{(0)}(\Gamma)\leq Cq^{k}
\end{equation}
for some $q, 0<q<1$. Here $J_k$ is the Bessel function of order $k$. Moreover, as observed in \cite{J},
$$
J_k(kr)=C_0 \sqrt{k} (r^2-1)^{-\frac{1}{4}} cos (-\frac{\pi}{4}+k(\sqrt{r^2-1}-arc cos \frac{1}{r}))+o(k^{-\frac{1}{2}})
$$
where $C_0$ is some constant and $lim k^{\frac{1}{2}} o(k^{-\frac{1}{2}})=0$ as $k\rightarrow 0$ uniformly on $[1.5, 2]$. Hence
$$
\int_{1.5}^2 J^2_k(kr)=C^2_0 k^{-1}\int_{1.5}^2 (r^2-1)^{-\frac{1}{2}} cos^2 (-\frac{\pi}{4}+k(\sqrt{r^2-1}-arc cos \frac{1}{r})r dr+o(k^{-1})=
$$
$$
C^2_0 k^{-1}\int_{\rho(1.5)}^{\rho(2)}  cos^2 (-\frac{\pi}{4}+k \rho)r (r^2-1)^{-1} d\rho +o(k^{-1})
$$
where $\rho(r)=\sqrt{r^2-1}- arc cos \frac{1}{r}$. It is easy to see that $\frac{ d\rho}{dt}=\frac{\sqrt{r^2-1}}{r}$. Therefore,
$$
\int_{1.5}^2 J^2_k(kr)\geq
C^{-1} k^{-1}\int_{\rho(1.5)}^{\rho(2)}  cos^2 (-\frac{\pi}{4}+k \rho) d\rho=
$$
$$
(Ck)^{-1}(\rho(2)-\rho(1.5)-
(2k)^{-1}(sin(-\frac{\pi}{2}+
2k\rho(2))-sin(-\frac{\pi}{2}+
2k\rho(1.5)))\geq (Ck)^{-1}
$$ 
for large $k$. Now from the definition of $u( ;k)$ we conclude that
\begin{equation}
\label{J2}
\|u( ;k)\|_{(0)}(\Omega_0) \geq C^{-1}k^{-\frac{7}{6}},
\end{equation}
where $\Omega_0=\{x: 1.5<|x|<2\}$. \eqref{J1} and \eqref{J2} suggest that the bound \eqref{stability0E} is nearly sharp (indeed, with the choice $u=u( ;k)$ the  difference of powers of $k$ on the left and right hand sides
is $\frac{1}{6}$). It is feasible that precise (but fairly complicated) representations of $J_k(kr)$ near $r=1$ (\cite{W}, p.252) will imply that the bound
\eqref{stability0E} is sharp.

Now we comment on the further content of the paper. 

In section 2 we derive Theorems 1.1, 1.2 by mapping $\Omega$ onto special domains where one can make use of standard energy estimates for second order equations. Similar ideas were used in
\cite{I2}, \cite {I3}, \cite{I5}. Since our equations are of the elliptic type these estimates have (high order) terms which are not positive. To handle these terms we split $u$ into "low frequency" and "high frequency" parts. "Low frequency" parts are dominated by positive terms containing
$k^2$, while "high frequency" parts remain in the right hand sides of
\eqref{stability0}, \eqref{stability1},\eqref{stability0om},\eqref{stability1om} and can be viewed as a priori constraints. The special domains are diffeomorphic to the unit ball and are used to introduce semi norms $|||u|||_{(1,k)}$ and ("low frequency") subspaces $H^{(1)}(\Omega;k)$.
In less precise forms of stability estimates
\eqref{stability0}, \eqref{stability0om} these semi norms and subspaces are not needed and estimates can be obtained in general domains
(not diffeomorphic to the unit ball) by the conventional partitioning into special domains.
In section 3 we consider annular domains where a construction of a "low frequency" part is explicit by using spherical harmonics. Otherwise the arguments are very similar to section 2, we only have to combine them with 
some basic differentiation and 
 harmonic analysis on the unit sphere. To bound "high frequency" part in
 \eqref{stability0}-
\eqref{stability1om} we transform the elliptic equation with the parameter $k$ into a hyperbolic equation and derive logarithmic stability bounds in the lateral (non hyperbolic) Cauchy problem for this equation. To do so we use the
Fourier-Bros-Iagolnitzer (FBI) transform to replace the hyperbolic equation with an elliptic equation without large parameter $k$ as first proposed by Robbiano \cite{R1}, \cite{R2}.
 For the resulting elliptic equation in a standard domain we use known Carleman estimates with a  simple suitable weight function. The crucial step in the proof is the choice of 
parameters in the FBI transform and the Carleman weight function. 
Observe that the idea of converting an hyperbolic equation into an elliptic one (and back) was conceived  and used already by Hadamard in 1930s. 
Theorem 1.4 has similarities with results in \cite{R2}, however a crucial difference is that we do not assume any homogeneous boundary condition on 
$\partial\Omega\times (-T,T)$ as in \cite{R2} and derive a  conditional logarithmic stability stability in a sub domain $\Omega_0$ but not in $\Omega$. As far as we know the closest
 bound on $\|u\|_{(0)}(\Omega_0)\times(-T_0,T_0)$ is obtained by similar methods in \cite{BKL}, however with smaller power of $|ln \varepsilon|$ and when 
$v_0=v_1=0$. 
 Also we obtain bounds \eqref{stability2space} which are more suitable for applications to the Helmholtz type equations. In \cite{V} there are (different) logarithmic stability bounds in the whole $\Omega$ for the continuation problem under homogeneous boundary conditions on $\partial\Omega$.  Finally, in concluding section 5
 we comment on challenging questions and outline possible
 new directions and links to source and coefficients identification problems.

\section{Proofs
 under "high frequency"  energy a priori constraints}

In this section we will prove Theorems 1.1 and 1.2.

{\bf Proof of Theorem 1.1}

We will first prove \eqref{stability1}. Then 
$\bar\Omega$ is the $C^2$ diffeomorphic  image of some standard set $\Omega(1)=
\{x: h(|x'|)\leq x_n \leq 1 \}$, where $h$ is a function in $C({\mathbb R})$ and in $C^2([0,2])$ satisfying the conditions:
$h=0$ on $(0,1)$, $h'(1)=h''(1)=0$, $0<h'$ on $(1,2)$, $h=1$ on $(2,+\infty)$, under a $C^2(\bar\Omega(1))$-diffeomorphism $y(x)$, so that
$y(\partial\Omega(1)\cap\{0\leq x_n<1\})\subset \Gamma$. Since the form of the second order equation \eqref{PD} (in particular the ellipticity) is invariant under such diffeomorphisms it suffices to give a proof only for $\Omega=\Omega(1)$ and $\Gamma=\partial \Omega(1)\cap \{x_n<1\}$. 

Observe that
$$
a_{nn}(\partial_n^2 u \partial_n \bar u +\partial_n^2 \bar u \partial_n  u) e^{-\tau x_n}=
\partial_n(a_{nn}|\partial_n u|^2 e^{-\tau x_n})
+\tau a_{nn}
|\partial_n u|^2 e^{-\tau x_n}
-(\partial_n a_{nn})|\partial_n u|^2 e^{-\tau x_n},
$$
\begin{equation}
\label{aun}
a_{jn}(\partial_j\partial_n u \partial_n \bar u +
\partial_j\partial_n \bar u \partial_n u) e^{-\tau x_n}=
\partial_j(a_{jn}|\partial_n u|^2 e^{-\tau x_n})
-(\partial_j a_{jn})
|\partial_n u|^2 e^{-\tau x_n},j=1,...,n-1.
\end{equation}
Let $\Omega(\theta)=\Omega\cap
\{0<x_n<\theta\}, 
\Gamma(\theta)=
\Gamma\cap\{0<x_n<\theta \}$.
Integrating by parts with respect to $x_j$, we yield
$$
\int_{\Omega(\theta)}\sum_{j=1}^{n-1}a_{jm}\partial_j\partial_m u \partial_n \bar u e^{-\tau x_n}=
\int_{\Gamma(\theta)}\sum_{j=1}^{n-1}
a_{jm}\partial_m u\partial_n \bar u e^{-\tau x_n} \nu_j d\Gamma-
$$
\begin{equation}
\label{au1}
\int_{\Omega(\theta)}
\sum_{j=1}^{n-1} a_{jm}\partial_m u \partial_j\partial_n \bar u e^{-\tau x_n}
-\int_{\Omega(\theta)}\sum_{j=1}^{n-1}(\partial_j a_{jm})\partial_m u \partial_n \bar u e^{-\tau x_n}.
\end{equation}
We have
$$
\sum_{j,m=1}^{n-1} a_{jm}(
\partial_m u \partial_j\partial_n \bar u + 
\partial_m \bar u \partial_j\partial_n u)
e^{-\tau x_n}=
\sum_{j,m=1}^{n-1}
\partial_n( a_{jm}\partial_m u \partial_j \bar u e^{-\tau x_n})+
$$
\begin{equation}
\label{sym}
\tau  \sum_{j,m=1}^{n-1} a_{jm}\partial_m u \partial_j \bar u e^{-\tau x_n}-
\sum_{j,m=1}^{n-1}(\partial_n a_{jm})\partial_m u \partial_j \bar u e^{-\tau x_n},
\end{equation}
due to the symmetry of $a_{jm}$. 

To form a standard  energy integral we multiply the both sides of \eqref{PD}
by $\partial_n \bar u e^{-\tau x_n}$, add its complex conjugate,  and integrate by parts over 
$\Omega(\theta), 0<\theta < 1,$ with using \eqref{aun}, \eqref{au1}, \eqref{sym}, and the notation $B(\theta)=\{x':
|x'|<h^{-1}(\theta)\}$ to yield
$$
\int_{B(\theta)}
a_{nn}|\partial_n u|^2( ,\theta) e^{-\tau \theta }+
\int_{\Gamma(\theta)}
a_{nn}|\partial_n u|^2 e^{-\tau x_n}\nu_n d\Gamma 
+ {\tau}\int_{\Omega(\theta)}
a_{nn}|\partial_n u|^2
 e^{-\tau x_n}+
$$
$$
2\int_{\Gamma(\theta)}
\sum_{j=1}^{n-1}a_{jn}
|\partial_n u|^2 \nu_j e^{-\tau x_n} d\Gamma+
\int_{\Gamma(\theta)}\sum_{j,m=1}^{n-1}a_{jm}
(\partial_m u \partial_n \bar u
+ \partial_m \bar u 
\partial_n u) 
e^{-\tau x_n}\nu_j d\Gamma -
$$
$$
\int_{B(\theta)}
\sum_{j,m=1}^{n-1}a_{jm}
\partial_j u \partial_m \bar u( , \theta) e^{-\tau \theta}-
\int_{\Gamma(\theta)}\sum_{j,m=1}^{n-1}a_{jm}
\partial_j u \partial_m \bar u e^{-\tau x_n}\nu_n d\Gamma -
 \tau\int_{\Omega(\theta)}
 \sum_{j,m=1}^{n-1}a_{jm}
 \partial_j u \partial_m \bar u e^{-\tau x_n}+
$$
$$
k^2\int_{B(\theta)}|u|^2( ,\theta) e^{-\tau \theta }+
k^2\int_{\Gamma(\theta)}|u|^2 e^{-\tau x_n}\nu_n d\Gamma 
+
\tau k^2 \int_{\Omega(\theta)}
 |u|^2 e^{-\tau x_n}+...=
 $$
\begin{equation}
\label{energyu}
 \int_{\Omega(\theta)}(\partial_n \bar u  f + \partial_n u  \bar f) e^{-\tau x_n},
 \end{equation}
where $...$ denotes the sum of terms bounded by
$$
C \int_{\Omega}
 (\sum_{j=1}^{n}|\partial_j u|^2+k^2 |u|^2) e^{-\tau x_n}.
$$

To bound the fifth and seventh  integrals on the left hand side of \eqref{energyu} we need to split $u$ into "low" and "high" frequencies parts.
To do so we will extend $u$ from $\Omega$ onto ${\mathbb R}^{n-1}\times (0,1)$ and use the (partial) Fourier transform
  $\cal F$ with respect to $x'$.

Let $u^*(x',x_n)\in H^1({\mathbb R}^{n-1})$ be an extension of the function $
u(x',x_n)$ with respect to
$x'$ from $B(x_n)$ onto
${\mathbb R}^{n-1}$ which admits the bound
  \begin{equation}
\label{u*u}
\|u^*( ,x_n)\|_{(0)}({\mathbb R}^{n-1})\leq 
C_e \|u( ,x_n)\|_{(0)}(B(x_n)),
\;0<x_n<1.
\end{equation}
Since the radii of balls $B(x_n)$ are in $(1,2)$ the standard extension operators satisfy the bound \eqref{u*u}.
 We introduce low and high frequency projectors
\begin{equation}
\label{ul*}
u_l(x)={\cal F}^{-1} \chi_k {\cal F} u^*(x),\; u_h=u-u_l,
\end{equation}
where $\cal F$ is the (partial) Fourier transformation with respect to
$x'=(x_1,...,x_{n-1},0)$, $\chi_k(\xi')=1$ when $C_e E_0^2|\xi'|^2<(1-\varepsilon_1)k^2$ for some positive $\varepsilon_1$ and $\chi_k(\xi')=0$ otherwise.

We have
$$
  \sum_{j,m=1}^{n-1}a_{jm}( ,x_n)\partial_j u \partial_m \bar u( , x_n)=
 \sum_{j,m=1}^{n-1} a_{jm}( ,x_n)\partial_j (u_l+u_h) \partial_m ( \bar u_l+ \bar u_h)( , x_n)=
  $$
  $$
 \sum_{j,m=1}^{n-1}(a_{jm}( ,x_n)\partial_j u_l \partial_m \bar u_l( , x_n)+
  2 a_{jm}( ,x_n)\partial_j u_l \partial_m \bar u_h( , x_n)+
 a_{jm}( ,x_n)\partial_j u_h \partial_m \bar u_h( , x_n))\leq
 $$
 \begin{equation}
 \label{uhbound}
 \sum_{j,m=1}^{n-1}a_{jm}( ,x_n)\partial_j u_l \partial_m \bar u_l( , x_n)+
  C\delta\sum_{j=1}^{n-1}|\partial_j u_l|^2( ,x_n)+
  C\delta^{-1}\sum_{j=1}^{n-1}|\partial_j u_h|^2( ,x_n),
  \end{equation}
  where we used the elementary inequality $|AB|\leq \frac{\delta}{2}|A|^2+
  \frac{1}{2\delta}|B|^2$ with $A=\partial_j u_l, B=\partial_m u_h$  and $\delta \in (0,1)$ to be chosen later. 
 
According to the definition of $E_0$,
$$
-\int_{B(x_n)}\sum_{j,m=1}^{n-1}a_{jm}( ,x_n)
\partial_j u_l \partial_m \bar u_l( , x_n)\geq 
-\int_{B(x_n)}E_0^2 \sum_{j=1}^{n-1}|\partial_j u_l|^2( , x_n)
$$
$$
\geq 
-\int_{{\mathbb R}^{n-1}}E_0^2 \sum_{j=1}^{n-1}|\partial_j u_l|^2( , x_n)=
-\int_{{\mathbb R}^{n-1}}E_0^2 \sum_{j=1}^{n-1}\xi_j^2|{\cal F} u_l|^2( , x_n)\geq
-\int_{{\mathbb R}^{n-1}}k^2(1-\varepsilon_1)C_e^{-1}|{\cal F} u_l|^2( , x_n)=
$$
\begin{equation}
\label{uuk}
-(1-\varepsilon_1)k^2 C_e^{-1}\int_{{\mathbb R}^{n-1}} |u_l|^2( , x_n) \geq
-(1-\varepsilon_1)k^2 C_e^{-1}\int_{{\mathbb R}^{n-1}} |u^*|^2( , x_n)
\geq
-(1-\varepsilon_1)k^2\int_{B(x_n)} |u|^2( , x_n),
\end{equation}
where we utilized the Parseval's identity and used that ${\cal F} u_l(,\xi',x_n)=0$ when $- C_e E_0^2|\xi'|^2<-(1-\varepsilon_1)k^2$, due to \eqref{ul*}, and \eqref{u*u}. Similarly,
\begin{equation}
\label{ulk}
\int_{B(x_n)}
 \sum_{j=1}^{n-1}|\partial_j u_l|^2 ( ,x_n) \leq
  C k^2 \int_{B(x_n)} 
  |u|^2 ( ,x_n).
 \end{equation}
 Therefore, using \eqref{uhbound} and \eqref{uuk} we obtain
 $$
 -\int_{B(x_n)}\sum_{j,m=1}^{n-1}a_{jm}( ,x_n)
\partial_j u_l \partial_m \bar u_l( , x_n)\geq
$$
 \begin{equation}
 \label{au}
 -(1-\varepsilon_1)k^2\int_{B(x_n)}|u|^2 ( , x_n)
-C \delta k^2\int_{B(x_n)}|u|^2 ( , x_n)
-\frac{C}{\delta}\int_{B(x_n)} \sum_{j=1}^{n-1}|\partial_j u_h|^2 ( , x_n).
\end{equation}

Hence from \eqref{energyu} and \eqref{au} by using the inequalities
$2 bc\leq b^2+c^2$ and $\frac{1}{C}<a_{nn}$ (due to the ellipticity of $A$) we conclude that
$$
\frac{1}{C}\int_{B(\theta)}|\partial_n u|^2( , \theta) e^{-\tau \theta }
+ \frac{\tau}{C}\int_{\Omega(\theta)}|\partial_n u|^2 e^{-\tau x_n}+
$$
$$
(\varepsilon_1-C\delta) k^2\int_{B(\theta)} |u|^2( ,\theta) e^{-\tau \theta}+ 
\tau (\varepsilon_1-C\delta)k^2 \int_{\Omega(\theta)} |u|^2 e^{-\tau x_n}\leq
$$
$$
C(\int_{\Gamma}(|\nabla u|^2+k^2 |u|^2) d\Gamma +
\int_{\Omega}|f|^2 +
$$
$$
\int_{B(\theta)}\frac{1}{\delta}\sum_{j=1}^{n-1}|\partial_j u_h|^2( ,\theta) e^{-\tau \theta}+
\int_{\Omega(\theta)}(\frac{\tau}{\delta}\sum_{j=1}^{n-1}|\partial_j u_h|^2+
\sum_{j=1}^{n}
|\partial_j u|^2+k^2|u|^2) e^{-\tau x_n}).
$$
Let  $\delta=\frac{\varepsilon_1}{2C}$ and use that $u=u_l+u_h$, then we yield the inequality
$$
\int_{B(\theta)}|\partial_n u|^2( , \theta) e^{-\tau \theta}
+ \tau\int_{\Omega(\theta)}|\partial_n u|^2 e^{-\tau x_n}+
k^2\int_{B(\theta)} |u|^2( ,\theta) e^{-\tau\theta}+ 
\tau k^2 \int_{\Omega(\theta)} |u|^2 e^{-\tau x_n}\leq
$$
$$
C(\int_{\Gamma}(|\nabla u |^2+k^2 |u|) d\Gamma+
\int_{\Omega}|f|^2 +
$$
\begin{equation}
\label{energyu1}
\int_{B(\theta)}\sum_{j=1}^{n-1}|\partial_j u_h|( ,\theta)^2 e^{-\tau \theta}+
\int_{\Omega(\theta)}(\tau \sum_{j=1}^{n-1}|\partial_j u_h|^2+\sum_{j=1}^{n-1}|\partial_j u_l|^2+|\partial_n u|^2+k^2|u|^2) e^{-\tau x_n}).
\end{equation}
Choosing and fixing sufficiently large $\tau$ (depending on the same parameters as $C$)
to absorb the three last terms on the right hand side in \eqref{energyu1} (with the help of \eqref{ulk}) by the left hand side we obtain
$$
\int_{B(\theta)}|\partial_n u |^2( , \theta)
+ \int_{\Omega(\theta)}|\partial_n u|^2 +
k^2\int_{B(\theta)} |u|^2( ,\theta) + 
k^2 \int_{\Omega(\theta)} |u|^2 \leq
$$
\begin{equation}
\label{energyu2}
C(\int_{\Gamma}(|\nabla u|^2 +k^2 |u|^2) d\Gamma +
\int_{\Omega}|f|^2 +\int_{B(\theta)}\sum_{j=1}^{n-1}|\partial_j u_h|^2( ,\theta)+
\int_{\Omega}\sum_{j=1}^{n-1}|\partial_j u_h|^2).
\end{equation}

Integrating the inequality \eqref{energyu2} with respect $\theta$ over $(0,1)$ and
dropping the first two terms on the left side,  we yield
$$
k^2\| u\|^2_{(0)}(\Omega)
 \leq
 C(\|u_1\|_{(0)}^2(\Gamma)+\|u_0\|_{(1)}^2(\Gamma)+k^2\|u_0\|_{(0)}^2(\Gamma)+ \|f\|_{(0)}^2(\Omega)+
\sum_{j=1}^{n-1}\|\partial_j u_h\|_{(0)}^2(\Omega)).
$$
Letting
$$
|||u|||_{(1,k)}(\Omega)=
(\sum_{j=1}^{n-1}||\partial_j u_h||^2_{(0)}(\Omega))^{\frac{1}{2}}
$$
 and dividing by $k^2$ we obtain \eqref{stability1}.

The bound \eqref{stability0}
follows from \eqref{stability1} by a suitable partitioning of $\Omega$. We claim that there is a finite covering
$\Omega_1,...,
\Omega_J$ of $\Omega$ such that $\Omega_j$ is the $C^2
(\bar\Omega(1))$-diffeomorphic image of $\Omega(1)$ and 
$\partial\Omega_j\cap \Gamma$  is a non empty
open (in $\partial\Omega$) subset of $\partial \Omega$.

Indeed, if $x\in \partial\Omega$ there is $\Omega^x$ which is the
$C^2(\bar\Omega(1))$-diffeomorphic image of $\Omega(1)$  such that
$\partial\Omega^x\cap\Gamma$
is a non empty open (in $\partial\Omega$) set 
 and $\partial \Omega^x\cap\partial\Omega$ is also open set containing
 $x$. If $x\in \Omega$ there
 is $\Omega^x$ which is
 the
$C^2(\bar\Omega(1))$-diffeomorphic image of $\Omega(1)$  such that
$\partial\Omega^x\cap\Gamma$
is a non empty open (in $\partial\Omega$) set 
 and $x\in \Omega^x$.
 $\partial\Omega^x\cap\partial\Omega$ 
 form an open covering of compact set $\partial\Omega$, so there is a finite sub covering
 $\Omega_1,...,\Omega_{J(1)}$. Then 
 $\Omega\setminus(\Omega_1\cup...
 \Omega_{J(1)})$ is compact and hence it has a finite sub covering $\Omega_{J(1)+1},...,\Omega_J$ by some $\Omega^x$. 

Obviously,
$$
\|u\|_{(0)}(\Omega)^2\leq
\|u\|_{(0)}(\Omega_1)^2+...
\|u\|_{(0)}(\Omega_J)^2\leq
$$
$$
C(\|u_0\|^2_{(0)}(\Gamma)
+\frac{1}{k^2}(\|f\|^2_{(0)}(\Omega)+
\|u_0\|^2_{(1)}(\Gamma)+
\|u_1\|^2_{(0)}(\Gamma)+
\|u\|^2_{(1)}(\Omega))),
$$
since
$$
\|u\|_{(0)}(\Omega_j)^2\leq
C(\|u_0\|^2_{(0)}(\Gamma\cap\partial\Omega_j)+
+\frac{1}{k^2}(\|f\|^2_{(0)}(\Omega_j)+
\|u_0\|^2_{(1)}(\Gamma\cap\partial\Omega_j)+
\|u_1\|^2_{(0)}(\Gamma\cap\partial\Omega_j)+
\|u\|^2_{(1)}(\Omega_j)))\leq
$$
$$
C(\|u_0\|^2_{(0)}(\Gamma)+
+\frac{1}{k^2}(\|f\|^2_{(0)}(\Omega)+
\|u_0\|^2_{(1)}(\Gamma)+
\|u_1\|^2_{(0)}(\Gamma)+
\|u\|^2_{(1)}(\Omega)
$$
due to \eqref{stability1}.

The proof is complete.

{\bf Proof of Theorem 1.2}

As in the proof of Theorem 1.1, we first prove
\eqref{stability1om}. Observe, that then $\Omega$ is the image of $\Omega(2)\subset
\{x: 0< x_n<1 \}$  under a $C^2(\bar\Omega(2))$-diffeomorphism $y(x)$, so that
$y(\partial\Omega(2)\cap
\{x_n=1\})\subset (\partial\Omega \setminus \bar\omega)$ and the image of $\partial \Omega(2) \cap \bar \omega(2)$ where
$\omega(2)$ is an open subset of $\Omega(2)$ with
$\partial\Omega(2)\cap \{x_n<1\} \subset \bar \omega(2)$ and the image of 
$\partial\Omega(2)\cap \bar \omega(2)$ is inside of an open  (in $\partial\Omega$) subset of
$\bar\omega\cap\partial\Omega$.
One can choose $\Omega(2)$ to be invariant with respect to rotations around the $x_n$-axis.  Since the form of the second order equation \eqref{PD} (in particular the ellipticity) is invariant under such diffeomorphisms it suffices to give a proof only for $\Omega=\Omega(2)$ and $\omega= \omega(2)$.  

We comment on a possible choice of such a diffeomorphism. We can assume that $\Omega$ is the ball
centred at $(0,...,0,\frac{1+\delta}{2})$ of the radius $\frac{1-\delta}{2}$ and $\omega$ contains
$\Omega\cap B_1$ where $B_1$ is some ball centred at $(0,0,...,\delta)$. Choosing small $\delta$ we can achieve that the image of $\Omega$ under the inversion
$|x|^{-2}x$, a translation in the $x_n$-direction and some  scaling is the unit ball $B$ and the closure of the image of $\omega$ contains
$\partial B\cap\{x_n<0.5\}$.
Now we use the same notation
$x$ for variables after transformation and "flatten"
the part of $\partial B$ outside the image of $\omega$.
To do so we can use the map
$$
x^*(x)=(x', (1-\chi_n(x_n) \chi_1(|x'|)x_n+
\chi_n(x_n)\chi_1(|x'|) |x|
$$
where $x'=(x_1,...,x_{n-1})$, $\chi_n \in C^{\infty}(\mathbb{R}),  \chi_n(x_n)=1$
when $0.5<x_n<1$, $\chi_n(x_n)=0$ when
$x_n<0.3$ and $0\leq \chi_n'$.
To define the function
$\chi_1 $ we pick up positive
numbers $\delta_1, \delta_2,\delta_3, \frac{\sqrt{3}}{2}<\delta_1 < \delta_2 <\delta_3 < 1$ and let $\chi_1\in C^{\infty}(\mathbb{R}), \chi_1(r)=1$ when $ r<\delta_1$, $\chi_1(r)=0$ when $\delta_3<r$ and $0\leq \chi_1
\leq 1$. Then
$$
\frac {\partial x^*_n}{\partial x_n}=
1-\chi_1\chi_n +
\chi_1\chi_n |x|^{-1}x_n+\chi_1 \chi_n'(|x|-x_n).
$$
Considering the cases when
$0.5<x_n$ (then $\frac {\partial x^*_n}{\partial x_n}
1-\chi_1+ \chi_1|x|^{-1}x_n $), when $0.3\leq x_n\leq 0.5$ and when $x_n<0.3$ (then
$\frac {\partial x^*_n}{\partial x_n}=1 $), we conclude that $0<\frac {\partial x^*_n}{\partial x_n}$ on $\{x: |x|\leq 1\}$, and hence $x^*(x)$ is a $C^{\infty}$-diffeomorphism
of the closed unit ball transforming the part of $\partial B$ outside $\omega$ into the hyperplane $\{x_n^*=1\}$. Obviously we have the rotational invariance around the $x_n$-axis.

Let $\chi$ be a $C^2(\mathbb R^n)$ cut off function,
such that $\chi=1$ on $\Omega\setminus \omega$, $\chi=0$ on $\mathbb R^{n-1}\times (0,1) \setminus \Omega$. 
Let $u_*=\chi u$.  From \eqref{PD} by using the Leibniz formula we yield 
\begin{equation}
\label{Au*}
(\sum_{j,m=1}^n a_{jm}\partial_j\partial_m  +
\sum_{j=1}^n a_j\partial_j  +a-ika_0+k^2)u_*= f_*,
\end{equation}
where 
$$
f_*= \chi f+2\sum_{j,m=1}^n a_{jm}\partial_j\chi \partial_l u +
(\sum_{j,m=1}^n a_{jm}\partial_j\chi \partial_m \chi+ \sum_{j=1}^n a_j\partial_j \chi) u.
$$

To form an  energy integral  we multiply the both sides of \eqref{Au*}
by $\partial_n \bar u_* e^{-\tau x_n}$, add its complex conjugate,  and integrate by parts over 
$\{ 0< x_n<\theta \}$ with using \eqref{aun}, \eqref{au1}, \eqref{sym},  and as in the proof of Theorem 1.1 we yield
$$
\int_{\mathbb R^{n-1}}
a_{nn}|\partial_n u_*|^2( ,\theta) e^{-\tau \theta } 
+ {\tau}\int_{\mathbb R^{n-1}\times (0,\theta)}
a_{nn}|\partial_n u_*|^2
 e^{-\tau x_n}-
 $$
 $$
\int_{\mathbb R^{n-1}}
\sum_{j,m=1}^{n-1}a_{jm}
\partial_j u_* \partial_m \bar u_*( , \theta) e^{-\tau \theta}-
 \tau\int_{\mathbb R^{n-1}\times (0,\theta)}
 \sum_{j,m=1}^{n-1}a_{jm}
 \partial_j u_* \partial_m \bar u_* e^{-\tau x_n}+
$$
$$
k^2\int_{\mathbb R^{n-1}}|u_*|^2( ,\theta) e^{-\tau \theta }
+
\tau k^2 \int_{\mathbb R^{n-1}\times (0,\theta)}
 |u_*|^2 e^{-\tau x_n}+...=
 $$
\begin{equation}
\label{energyu*}
 \int_{\mathbb R^{n-1}\times (0,\theta)}(\partial_n \bar u_*  f_* + \partial_n u_*  \bar f_*) e^{-\tau x_n},
 \end{equation}
where $...$ denotes the sum of terms bounded by
$$
C \int_{{\mathbb R}^{n-1}\times(0,1)}
 (\sum_{j=1}^{n}|\partial_j u_*|^2+k^2|u_*|^2) e^{-\tau x_n}.
$$

As in the proof of Theorem 1.1 we  split $u_*$ into "low" and "high" frequencies parts.
To do so we will extend $u_*$ from $\Omega$ onto $({\mathbb R}^{n-1}\times (0,1))\setminus \Omega$ as zero and use the (partial) Fourier transform
  $\cal F$ with respect to $x'$.
 We introduce low and high frequency projectors
\begin{equation}
\label{ul**}
u_{*l}(x)={\cal F}^{-1} \chi_k {\cal F} u_*(x),\; u_{*h}=u_*-u_{*l},
\end{equation}
where $\chi_k(\xi')=1$ when $E_0^2|\xi'|^2<(1-\varepsilon_1)k^2$ for some positive $\varepsilon_1$ and $\chi_k(\xi')=0$ otherwise.

As in \eqref{uhbound},
$$
  \sum_{j,m=1}^{n-1}a_{jm}( ,x_n)\partial_j u_* \partial_m \bar u_*( , x_n)
  \leq
 $$
 \begin{equation}
 \label{uhbound*}
 \sum_{j,m=1}^{n-1}a_{jm}( ,x_n)\partial_j u_{*l} \partial_m \bar u_{*l}( , x_n)+
  C\delta\sum_{j=1}^{n-1}|\partial_j u_{*l}|^2( ,x_n)+
  C\delta^{-1}
  \sum_{j=1}^{n-1}|\partial_j u_{*h}|^2( ,x_n).
  \end{equation}
By using \eqref{ul**}, similarly to
\eqref{uuk}, \eqref{ulk}
we have
\begin{equation}
\label{uuk*}
-\int_{{\mathbb R}^{n-1}}\sum_{j,m=1}^{n-1}a_{jm}( ,x_n)
\partial_j u_{*l} \partial_m \bar u_{*l}( , x_n)\geq 
-(1-\varepsilon_1)k^2\int_{{\mathbb R}^{n-1}} |u_{*}|^2( , x_n),
\end{equation}
\begin{equation}
\label{ulk*}
\int_{{\mathbb R}^{n-1}}
 \sum_{j=1}^{n-1}|\partial_j u_{*l}|^2 \leq
  Ck^2 \int_{{\mathbb R}^{n-1}} |u_{*}|^2.
 \end{equation}
 From \eqref{uhbound*},  \eqref{uuk*} and \eqref{ulk*} we obtain
 $$
 -\int_{{\mathbb R}^{n-1}}\sum_{j,m=1}^{n-1}a_{jm}( ,x_n)
\partial_j u_{*l} 
\partial_m \bar u_{*l}( , x_n)\geq
$$
 \begin{equation}
 \label{au*}
 -(1-\varepsilon_1)k^2\int_{{\mathbb R}^{n-1}}|u_{*}|^2 ( , x_n)
-C \delta k^2\int_{{\mathbb R}^{n-1}}|u_{*}|^2 ( , x_n)
-\frac{C}{\delta}\int_{{\mathbb R}^{n-1}} \sum_{j=1}^{n-1}|\partial_j u_{*h}|^2 ( , x_n).
\end{equation}

From \eqref{energyu*} and \eqref{au*}  we conclude that
$$
\frac{1}{C}\int_{{\mathbb R}^{n-1}}|\partial_n u_*|^2( , \theta) e^{-\tau \theta }
+ \frac{\tau}{C}\int_{{\mathbb R}^{n-1}\times(0,\theta)}|\partial_n u_{*}|^2 e^{-\tau x_n}+
$$
$$
(\varepsilon_1-C\delta) k^2\int_{{\mathbb R}^{n-1}} |u_{*}|^2( ,\theta) e^{-\tau \theta}+ 
\tau (\varepsilon_1-C\delta)k^2\int_{{\mathbb R}^{n-1}\times(0,\theta)} |u_{*}|^2 e^{-\tau x_n}\leq
$$
$$
C\int_{\Omega}|f_*|^2 e^{-\tau x_n}+
\int_{{\mathbb R}^{n-1}}\frac{1}{\delta}\sum_{j=1}^{n-1}|\partial_j u_{*h}|^2( ,\theta) e^{-\tau }+
\int_{{\mathbb R}^{n-1}\times(0,\theta)}(\frac{1}{\delta}\sum_{j=1}^{n-1}|\partial_j u_*|^2+
|\partial_n u_*|^2+k^2|u_*|^2) e^{-\tau x_n}).
$$
Let  $\delta=\frac{\varepsilon_1}{2C}$, then we yield the inequality
$$
\int_{{\mathbb R}^{n-1}}|\partial_n u_{*}|^2( , \theta) e^{-\tau \theta}
+ \tau\int_{{\mathbb R}^{n-1}\times(0,\theta)}|\partial_n u_*|^2 e^{-\tau x_n}+
k^2\int_{{\mathbb R}^{n-1}} |u_*|^2( ,\theta) e^{-\tau\theta}+ 
\tau k^2 \int_{{\mathbb R}^{n-1}\times(0,1)} |u^*|^2 e^{-\tau x_n}\leq
$$
\begin{equation}
\label{energyu*1}
C(
\int_{\Omega}|f_*|^2 e^{-\tau x_n}+
\int_{{\mathbb R}^{n-1}}\sum_{j=1}^{n-1}|\partial_j u_{*h}|^2( ,\theta) e^{-\tau}+
\int_{{\mathbb R}^{n-1}\times(0,1)}(\sum_{j=1}^{n-1}|\partial_j u_{*h}|^2+|\partial_n u_*|^2+k^2|u^*|^2) e^{-\tau x_n}).
\end{equation}
Choosing and fixing sufficiently large $\tau$ (depending on the same parameters as $C$)
to absorb the two last terms on the right hand side in \eqref{energyu*1} by the left hand side we obtain
$$
\int_{{\mathbb R}^{n-1}}|\partial_n u_*|^2( , \theta)
+ \int_{{\mathbb R}^{n-1}\times(0,\theta)}|\partial_n u_*|^2 +
k^2\int_{{\mathbb R}^{n-1}} |u_*|^2( ,\theta) + 
k^2 \int_{{\mathbb R}^{n-1}\times(0,\theta)} |u_*|^2 \leq
$$
\begin{equation}
\label{energyu*2}
C(\int_{\Omega}|f_*|^2 +\int_{{\mathbb R}^{n-1}}\sum_{j=1}^{n-1}|\partial_j u_{*h}|^2( ,\theta))+
\int_{{\mathbb R}^{n-1}\times(0,1)}\sum_{j=1}^{n-1}|\partial_j u_{*h}|^2).
\end{equation}

Integrating the inequality \eqref{energyu*2} with respect to $\theta$ over $(0,1)$, dropping the first two terms on the left hand side, and recalling that $u_*=\chi u$ we yield
$$
k^2\| u\|^2_{(0)}(\Omega)
 \leq
 C( \|f\|_{(0)}^2(\Omega)+ 
\|u\|_{(1)}^2(\omega)+
\sum_{j=1}^{n-1}\|\partial_j u_{*h}\|_{(0)}^2({\mathbb R}^{n-1}\times(0,1)).
$$
Letting
$$
|||u|||_{(1,k)}(\Omega)=
\sum_{j=1}^{n-1}\|\partial_j u_{*h}\|_{(0)}^2({\mathbb R}^{n-1}\times(0,1))
$$

using that $u=v$ on $\Omega\setminus V$, and dividing by $k^2$ we obtain \eqref{stability1om}.

\eqref{stability0om} follows from \eqref{stability1om} by partitioning $\Omega$ into the union of $C^2$ diffeomorphic images of
$\Omega(2)$ as in the proof of Theorem 1.1

\section{Proof for annular domains}

In this section we will prove Theorem 1.3. We will use polar coordinates
$(r,\sigma), \sigma\in S^n,$
and the operator $A$ in these coordinates:
\begin{equation}
 \label{Apolar}
Au=a^{nn}\partial_r^2 u+A_{1,\sigma}\partial_r u+
A_{2,\sigma} u
\end{equation}
where $A_{j,\sigma}$ is a $j$-th order linear partial
differential operator on $S^n$.

Let $\nabla_{\sigma}$ be the tangential gradient on $S^n$, i.e. the orthogonal projection of the gradient onto the tangent space and 
the normal gradient $\nabla_{nu}=\nabla-\nabla_{\sigma}$.
We also will use the tangential and normal divergences of the vector field $V=(V_1,...,V_n)$ defined as
$$
\nabla_\sigma\cdot V=
\sum_{j=1}^n\nabla_{\sigma}
V_j\cdot e(j),\;\nabla_\nu\cdot V=
\sum_{j=1}^n\nabla_{\nu}
V_j\cdot e(j),
$$
where $e(1),...,e(n)$ is the standard orthonormal base in $\mathbb R^n$.
By using local coordinates one can see that  
$$
a^{nn}(\partial_r)^2u=div_{\nu}(a \nabla_{\nu} u) +..., A_{2,\sigma}u =\nabla_{\sigma} \cdot (a \nabla_{\sigma} u)+...,
$$
 where $...$ are terms with 
the absolute value bounded by $C(|\nabla_{\sigma} u|+|u|)$.
We recall the integration by parts formula on $S^n$:
\begin{equation}
 \label{stokes}
\int_{S^n} g \nabla_{\sigma} \cdot V dS=-\int_{S^n} V\cdot \nabla_{\sigma} g dS
\end{equation}
for a vector field $V$ and a function $g$ on $S^n$.

From \eqref{PD}, \eqref{Apolar} we yield
\begin{equation}
\label{PDEup}
a^{nn}\partial_r^2 u+A_{1,\sigma}\partial_r u+
A_{2,\sigma} u+ika_0u+k^2u=
 f \;\mbox{in}\; S^n\times(1,R).
\end{equation}

By using \eqref{stokes} we have
$$
\int_{S^n} (A_{2,\sigma} u\partial_r\bar u +A_{2,\sigma} \bar u\partial_r u)dS=
\int_{S^n} ((\nabla_{\sigma} \cdot (a \nabla_{\sigma} u)) \partial_r\bar u + (\nabla_{\sigma}\cdot (a \nabla_{\sigma} \bar u)) \partial_r  u)+...)dS=
$$
\begin{equation}
\label{A2u}
-\int_{S^n} ((a\nabla_{\sigma} u)\cdot  \nabla_{\sigma}\partial_r\bar u + 
(a\nabla_{\sigma}\bar u)\cdot  \nabla_{\sigma}\partial_r u+...)dS=
-\int_{S^n} (\partial_r(a\nabla_{\sigma} \cdot \nabla_{\sigma}\bar u)+...)dS,
\end{equation}

Repeating the argument from the proof of Theorem 1.1 ( multiplying the both parts of \eqref{PDEup} by $\partial_r \bar u e^{-\tau r}$, adding its complex conjugate, and integrating by parts over 
$S^n\times(1,\rho), 1<\rho\leq R$), we will have
$$
\int_{S^n}a^{nn}|\partial_r u|^2( ,\rho) e^{- \tau \rho} \rho^{n-1} dS-
\int_{S^n}a^{nn}|\partial_r u|^2( ,1) e^{-\tau}dS 
+ {\tau}\int_{S^n\times(1,\rho)}a^{nn}|\partial_r u|^2 e^{-\tau r}r^{n-1} dr dS-
$$
$$
\int_{S^n}a \nabla_{\sigma} u \cdot  \nabla_{\sigma}\bar u( ,\rho)  e^{-\tau \rho}\rho^{n-1}dS+
\int_{S^n} a \nabla_{\sigma} u \cdot  \nabla_{\sigma}\bar u( ,1) e^{ -\tau}dS-
 \tau\int_{S^n\times(1,\rho)} a \nabla_{\sigma} u \cdot  \nabla_{\sigma}\bar u  e^{-\tau r}r^{n-1} dr dS+
$$
$$
k^2\int_{S^n}|u|^2( , \rho) e^{-\tau \rho}\rho^{n-1} dS-
k^2\int_{S^n}|u|^2( , 1) dS 
+
\tau k^2\int_{S^n\times(1,\rho)}|u|^2 e^{-\tau r}r^{n-1}dr dS+...=
$$
\begin{equation}
\label{energyup}
 \int_{S^n\times(1,\rho)}(\partial_r \bar u f+\partial_r u \bar f) e^{-\tau r}r^{n-1} dr dS,
\end{equation}
where $...$ denotes the sum of terms bounded by
$$
C \int_{\Omega}(|\nabla u|^2+k^2 |u|^2) e^{-\tau r}.
$$

To handle the fourth and sixth terms on the left hand side of \eqref{energyup} we use that
$$
-\int_{S^n}a\nabla_{\sigma}u
\cdot  \nabla_{\sigma}\bar u ( ,r) dS\geq
-\int_{S^n}a
\nabla_{\sigma}u_l
\cdot  \nabla_{\sigma}\bar u_l( ,r)dS
-\delta\int_{S^2}
|\nabla_{\sigma}u_l|^2 ( ,r)
dS -\frac{C}{\delta}
\int_{S^2}|\nabla_{\sigma}u_h|^2( ,r)dS.
$$

 As in the proof of Theorem 1.1, using \eqref{ulp} and that $e(\sigma;m,p)$ are eigenfunctions of the Beltrami operator  we yield
$$
-\int_{S^n}a\nabla_{\sigma}
u_l \cdot  \nabla_{\sigma}\bar u_l ( ,r) dS\geq
-E_0^2\int_{S^n} |\nabla_{\sigma} u_l|^2 ( ,r)dS=
$$
$$
r^{-2}E_0^2\int_{S^n} \Delta_{\sigma} u_l \bar u_l ( ,r)dS\geq 
-(1-\varepsilon)r^{-2}
\int_{S^n} |u|^2( ,r)dS
\geq 
-(1-\varepsilon)
\int_{S^n} |u|^2( ,r)dS,
$$
when $1<r$ and similarly
\begin{equation}
\label{boundulp}
-\int_{S^n} |\nabla_{\sigma} u_l|^2 ( ,r)\geq
-Ck^2\int_{S^n}|u|^2 ( ,r)dS,\;1<r.
\end{equation}

Hence from \eqref{energyup} we conclude that
$$
\int_{S^n}a^{nn}|\partial_r u|^2( , \rho) e^{-\tau \rho } \rho^{n-1} dS
+ \tau\int_{S^n\times(1,\rho)}
a^{nn}|\partial_r u|^2 e^{-\tau r}r^{n-1}dr dS+
$$
$$
(\varepsilon-C\delta) k^2 \int_{S^n} |u|^2( ,\rho) e^{-\tau \rho} \rho^{n-1}+ 
\tau (\varepsilon-C\delta) k^2
\int_{S^n\times(1,\rho)} |u|^2 e^{-\tau r}r^{n-1} dr dS\leq
$$
$$
C(\int_{S^n}(|\partial_r u|^2( ,1)+k^2 |u|^2( ,1)) dS+
\int_{S^n\times(1,R)}|f|^2 e^{-\tau r}r^{n-1} dr dS+
$$
$$
 \frac{C}{\delta}(\int_{S^n}(|\nabla_{\sigma} u_h|^2)( ,\rho) e^{-\tau \rho}\rho^{n-1}+
\tau \int_{S^n\times(1,\rho)}(|\partial_r u_h|^2+|\nabla_{\sigma} u_h|^2) e^{-\tau r}r^{n-1})dr dS+
 $$
 $$
\int_{S^n\times(1,\rho)}(|\nabla u|^2 +
 k^2|u|^2) e^{-\tau r}r^{n-1}dr dS).
$$
Choosing $\delta=\frac{\varepsilon}{2C}$ and using that $\frac{1}{C}<a^{nn}$ we yield
$$
\int_{S^n}(|\partial_r u|^2( , \rho)+k^2 |u|^2( ,\rho)) e^{-\tau \rho } \rho^{n-1} dS
+ \tau\int_{S^n\times(1,\rho)}(|\partial_r u|^2 +k^2 |u|^2) e^{-\tau r}r^{n-1} dr dS\leq
$$
$$
C(\int_{S^n}(|\partial_r u|^2( ,1)+k^2 |u|^2( ,1))dS+
\int_{S^n\times(1,R)}|f|^2 e^{-\tau r}r^{n-1}dr dS+
$$
\begin{equation}
\label{energyu1p}
\int_{S^n}|\nabla_{\sigma}u_h|^2( ,\rho) e^{-\tau \rho}\rho^{n-1}dS+
\tau \int_{S^n\times(1,\rho)}
|\nabla_{\sigma} u_h|^2 e^{-\tau r}r^{n-1}dr dS+
\int_{S^n\times(1,\rho)}(|\nabla u|^2 +
 k^2|u|^2) e^{-\tau r}r^{n-1} dr dS).
\end{equation}
Since $u=u_l+u_h$, from the triangle inequality we have
$$
\int_{S^n\times(1,\rho)}|\nabla_{\sigma} u|^2e^{-\tau r}r^{n-1}dr dS\leq 
2\int_{S^n\times(1,\rho)}(|\nabla_{\sigma} u_l|^2+|\nabla_{\sigma} u_h|^2)
e^{-\tau r}r^{n-1}dr dS\leq
$$
$$
 Ck^2 \int_{S^n\times(1,\rho)} |u|^2+
\int_{S^n\times(1,\rho)}|\nabla_{\sigma} u_h|^2 e^{-\tau r}r^{n-1}dr dS,
 $$
 when we apply \eqref{boundulp}. So from \eqref{energyu1p} we obtain
 $$
\int_{S^n}(|\partial_r u|^2( , \rho)+k^2 |u|^2( ,\rho)) e^{-\tau \rho } \rho^{n-1}dS
+ \tau\int_{S^n\times(1,\rho)}(|\partial_r u|^2 +k^2 |u|^2) e^{-\tau r}r^{n-1} dr dS\leq
$$
$$
C(\int_{S^n}(|\partial_r u|^2( ,1)+k^2 |u|^2( ,1))dS+
\int_{S^n\times(1,R)}|f|^2 e^{-\tau r}r^{n-1}drdS+
$$
$$
\int_{S^n}(|\nabla_{\sigma} u_h|^2)( ,\rho) e^{-\tau \rho}\rho^{n-1}dS+
\tau \int_{S^n\times(1,\rho)}(
|\nabla_{\sigma} u_h|^2)e^{-\tau r}r^{n-1}dr dS+
$$
$$
\int_{S^n\times(1,\rho)}(|\partial_r u|^2 + k^2|u|^2) e^{-\tau r}r^{n-1}dr dS).
$$
Now, choosing and fixing  $\tau$ sufficiently large (but depending on the same quantities as $C$) 
to absorb the  last term on the right side by the left side we yield
$$
\int_{S^n}(|\partial_r u|^2( , \rho)+k^2 |u|^2( ,\rho))dS
+ \int_{S^n\times(1,\rho)}(|\partial_r u|^2 +k^2 |u|^2)dr dS\leq
$$
$$
C(\int_{S^n}(|\partial_r u|^2( ,1)+k^2 |u|^2( ,1))dS+
\int_{\Omega}|f|^2 )+
$$
\begin{equation}
\label{energyu2p}
\int_{S^n}|\nabla_{\sigma} u_h|^2( ,\rho) dS +\int_{S^n\times(1,\rho)}|\nabla_{\sigma} u_h|^2 dr dS).
\end{equation}

By trace theorems for Sobolev spaces
$$
\|u_h( ,R)\|_{(1)}(\Gamma_1)\leq C(\theta)\|u_h\|_{(\frac{3}{2}+\theta)}(\Omega).
$$
For the high frequency part
\begin{equation}
\label{uhkp}
\|u_h\|_{(\frac{3}{2}+\theta)}^2(\Omega)
\leq Ck^{2\theta-1}\|u_h\|_{(2)}^2(\Omega),
\end{equation}
so from \eqref{energyu2p} we obtain \eqref{stability3}.

Integrating \eqref{energyu2p} with respect to $\rho$ over
$(1,R)$ we obtain
\begin{equation}
\label{energyu2pOmega}
 \int_{\Omega}(|\nabla u|^2 +k^2 |u|^2)\leq
C(\int_{\Gamma}(|\partial_r u|^2+k^2 |u|^2)d\Gamma+
\int_{\Omega}|f|^2 +
 \int_{\Omega}|\nabla_{\sigma} u_h|^2 )
\end{equation}
and hence derive \eqref{stability4}.

To obtain \eqref{stability5} we simply divide the both parts in \eqref{energyu2pOmega} by $k^2$.

The proof is complete.

\section{Proof for hyperbolic equations}

In this section we will prove Theorem 1.4.

By using compactness of $\bar\Omega_0$ we can cover this set by finitely many $\Omega_{0j}$, which are
 the images of an open set $\Omega^*_0=
\{\frac{1}{2}<y_n<1, |y|<1\}$ under $C^2$-diffeomorphisms $x( ;j)$ of the semi ball $B^+=\{y:|y|\leq 1, 0\leq y_n\}$ with $\overline{x( ;j)B^+}\subset
\Omega\cup\Gamma\cup(V\cap \partial\Omega)$. It suffices to consider two cases:
a) $ x( ;j)(\omega_0)\subset \omega $, where $\omega_0= \{y: 1-\delta<|y|<1\} \cap\{0<y_n\}$
with some positive $\delta$ and b) $x( ;j)(\overline{\partial B^+\cap\{0<y_n<1\}})\subset \Gamma$.
For brevity we return to the old notation replacing $y$ by $x$, so it is sufficient
to prove Theorem 1.3 when $\Omega$ is $B^+$, $\omega$ is $\omega_0$ or
$\Gamma$ is $\partial B^+\cap \{0<x_n\leq 1\}$. 

First we handle the case a).

Let $\chi_0$ be a $C^{\infty}({\mathbb R}^{n+1})$ function, $\chi_0=1$ on $ (\Omega\setminus\omega)\times
{\mathbb R}$, $\chi_0=0$ on $(\{0<x_n<1\}\setminus \Omega)\times {\mathbb R}$,
$0\leq \chi_0\leq 1$, and $\partial_t \chi_0=0$. Introducing
\begin{equation}
\label{w}
w=\chi_0 v,
\end{equation}
from (\ref{PDH}) we will have
\begin{equation}
\label{PDEw}
(A-a_0\partial_t -\partial_t^2) w= f_0,\;f_0=\chi_0 f+
2\sum_{j,m=1}^n a_{jm}\partial_j\chi_0\partial_m v+(A\chi_0-a \chi_0)v\;\mbox{in}\; Q=\Omega\times(-T,T).
\end{equation}

Let $\chi_{T}(t)$ be a $C^{\infty}({\mathbb R})$ function, $\chi_{T}(t)=1$ when $|t|<T-1$,
$0\leq\chi_{T}(t)\leq 1$, and $\chi_{T}(t)=0$ when $T<|t|<+\infty$.
As first suggested in \cite{R1}, we will make use of the FBI transform 
\begin{equation}
\label{FBI}
W(x,s;\lambda,t_1)=
\sqrt{\frac{\lambda}{2\pi}}\int_{{\mathbb R}}e^{-\frac{\lambda}{2}(is+t_1-t)^2}\chi_{T}(t)w(x,t) dt
\end{equation}
of a function $w(x,t)$.

Let $Q^s=\Omega\times(-1,1)$ be the domain in the $(x,s)$-space.
Using the time independence of the coefficients of $A$, applying the FBI transform to the both sides of  \eqref{PDEw}  and integrating by parts we obtain
\begin{equation}
\label{AW}
AW+a_0i \partial_s W+\partial_s^2 W= F_0 -F( ;T)\;\textit{on}\;Q^s,
\end{equation}
where
$$
 F(x,s ;T)=
\sqrt{\frac{\lambda}{2\pi}}\int_{{\mathbb R}}e^{-\frac{\lambda}{2}(is+t_1-t)^2}(a_0\partial_t\chi_{T}(t)w(x,t)
+2\partial_t\chi_{T}(t)\partial_t w(x,t)+\partial^2_t\chi_{T}(t)w(x,t))dt.
$$
Observe that, as in \cite{RZ},
$$
\|W\|^2_{(0)}(Q^s) \leq
C\lambda e^{\lambda}\|w||^2_{(0)}(Q), 
$$
 \begin{equation}
\label{Uu}
 \|\partial_jW\|_{(0)}^2(Q^s)
 \leq
C \lambda e^{\lambda}( \|\partial_j w\|_{(0)}^2(Q)+(1+\lambda)\| w\|_{(0)}^2(Q)), j=0,1,...,n,
\end{equation}
$$
 \|\partial_sW\|_{(0)}^2(Q^s)
 \leq
C \lambda^3 e^{\lambda}
 \| w\|_{(0)}^2(Q),
$$
$$
\|F( ;T)\|_{(0)}(Q^s)\leq
C\lambda^3 e^{\lambda-\lambda(T-1-T_1)^2} \| w\|^2_{(0)}(Q),
$$
provided $|t_1|<T_1<T-1$, where $T_1$ is to be chosen later. Here we let $\partial_0=\partial_s$ in the $(x,s)$ space and $\partial_0=\partial_t$ in 
the  $(x,t)$-space. Observe that to get the last bound \eqref{Uu} we integrated by parts in the integral defining $F( ;T)$ to eliminate $\partial_t w$. 
Similar bounds hold when we replace $\Omega$ by $\Gamma$.

We will use the Carleman weight
function
\begin{equation}
\label{weight}
\varphi(x,s)=e^{\gamma(x_n-s^2)}-1.
\end{equation}
As known \cite{I1}, section 3.2, there is $\gamma$ (depending only on $A, a_0, \Omega$) such that the following Carleman type estimate holds
\begin{equation}
\label{Carleman}
\int_{Q^s} \tau^{3-2|\alpha|}|\partial^{\alpha}W|^2 e^{2\tau\varphi} \leq
C(\int_{Q^s}|(A+a_0 i\partial_s+\partial_s^2)W|^2 e^{2\tau\varphi}+
\int_{\partial Q^s}(\tau^3 |W|^2+\tau|\nabla_{x,s} W|^2) e^{2\tau\varphi}),\;|\alpha|\leq 1,\:C<\tau.
\end{equation}

We will denote $Q^s(\delta)=
Q^s\cap\{\delta<\varphi\}$ and use another cut off function
$\chi( ;\delta)=0$ on $Q^s(\frac{\delta}{2})$ and $\chi( ;\delta)=1$ on $Q^s(\delta)$ with $|\partial^{\alpha}\chi( ;\delta)|\leq C \delta^{-|\alpha|}, |\alpha|\leq 2$.
Observe that $\chi( ;\delta)=0$  near $\Omega\times\{-1, 1\}$ and $\{x_n=0\}$.

Since
$$
(A+a_0 i\partial_s+\partial_s^2)(\chi W)=
$$
$$
\chi(A+a_0 i\partial_s+\partial_s^2) W+
2\partial_s \chi \partial_sW+
2\sum_{j,k=1}^n a^{jk}\partial_j \chi \partial_k W+((A-a+a_0i \partial_s +\partial_s^2)\chi) W,
$$
using $\chi( ;\delta_1)W$ instead of $W$ in \eqref{Carleman}, \eqref{PDEw}, and \eqref{AW} we yield
$$
\int_{Q^s(2\delta_1)} \tau^{3-2|\alpha|}|\partial^{\alpha}W|^2 e^{2\tau\varphi} \leq
C(\int_{Q^s}|F|^2 e^{2\tau\varphi}+\int_{ \omega^s}(|V|^2+|\nabla_x V|^2) e^{2\tau\varphi}+
$$
$$
\int_{Q^s}|F( ;T)|^2 e^{2\tau\varphi}+\int_{ Q^s\setminus Q^s(\delta_1)}(|W|^2+|\nabla_{x,s} W|^2) e^{2\tau\varphi}),\; |\alpha|\leq 1.
$$
Here $ \delta_1$ is a positive number to be chosen later, after \eqref{Wv}.

Let $\Phi=sup\varphi$ over $Q^s$. Since $\varphi\leq \Phi$ on $Q^s$ and $\delta <\varphi $ on $Q^s(\delta)$, it follows that
$$
e^{4\tau\delta_1}
\int_{Q^s(2\delta_1)} \tau^{3-2|\alpha|}|\partial^{\alpha}W|^2  \leq
C(e^{2\tau\Phi}(\|F^2\|_{(0)}(Q^s) + \| V\|^2_{(0)}(\omega^s))+ \| \nabla_x V\|^2_{(0)}(\omega^s))+
$$
$$
e^{2\tau\Phi}\|F( ;T)\|_{(0)}^2 (Q^s)+
e^{2\tau \delta_1}\|W\|_{(1)}^2(Q^s)).
$$
Using the bounds \eqref{Uu}  and recalling that 
$$
\varepsilon=\|f\|_{(0)}(Q)+
\|v_0\|_{L^2((-T,T);H^1(\Gamma))}+ \|v_1\|_{(0)}(\Gamma_0\times(-T, T))
+\|v\|_{L^2((-T,T);H^1(\omega))}
$$ 
we obtain
$$
\int_{Q^s(2\delta_1)}
|\partial^{\alpha}W( ;\lambda,t_2)|^2  \leq
C(e^{2\tau(\Phi+1)}e^{2\lambda}
\varepsilon^2+
$$
 \begin{equation}
\label{stabU}
e^{2\tau\Phi} e^{2\lambda}
e^{-\lambda(T-1-T_1)^2}\|v\|_{(0)}^2 (Q)+
e^{-2\tau \delta_1}e^{2\lambda}\|v\|_{L^2((-T,T);H^1(\Omega))}),
\end{equation}
provided $|t_2|<T_1<T-1$.

From the mean value bounds for the complex analytic function
$(\partial^{\alpha} W(x,s;\lambda,t_2))^2$ of
$z=t_2+is$ we have
$$
|\partial^{\alpha} W(x,0;\lambda,t_1)|^2\leq \frac{4}{\pi }\int_{\{|(t_2-t_1)+is|<\frac{1}{2}\}}|\partial^{\alpha} W(x,s;\lambda,t_2)|^2 ds dt_2\leq
$$
\begin{equation}
\label{stabU0}
\frac{4}{\pi}\int_{\{|t_2-t_1|<\frac{1}{2}, |s|<\frac{1}{2}\}}|\partial^{\alpha} W(x,s;\lambda,t_2)|^2 ds dt_2,
\end{equation}
when we assume that
$|t_1|<T_1-\frac{1}{2}$.

Now we make a crucial choice of 
\begin{equation}
\label{taulambda}
\tau=-\mu ln \varepsilon,\; \lambda = -\beta \mu ln \varepsilon, \varepsilon<1,
\end{equation}
where positive $\beta, \mu,
T, T_1, \kappa_1$
(depending only on $\Omega,T_0, A, a_0$) are selected so that
\begin{equation}
\label{betakappa}
\kappa_1+\beta<\delta_1,\;
\kappa_1+\mu((\Phi+1)+\beta)<1,\;
2(\kappa_1+\Phi)<\beta((T-1-T_1)^2-2).
\end{equation}
Due to this choice,
$$
e^{2\tau(\Phi+1)+2\lambda}
\varepsilon^2=
e^{-2\mu(\Phi+1) ln \varepsilon -2\beta\mu ln \varepsilon+2 ln \varepsilon}\leq e^{2\kappa_1 ln\varepsilon}=
\varepsilon^{2\kappa_1},
$$
$$
e^{-2\tau\delta_1+2\lambda}=
e^{2\mu(\delta_1-\beta) ln \varepsilon}\leq
e^{2\mu \kappa_1 ln\varepsilon}=
\varepsilon^{2\kappa_1\mu},
$$
$$
e^{2\tau \Phi+2\lambda-\lambda(T-1-T_1)^2}=
e^{\mu(2\Phi+2\beta-\beta(T-1-T_1)^2)(-ln \varepsilon)}\leq
e^{2\mu \kappa_1 ln\varepsilon}=
\varepsilon^{2\kappa_1\mu},
$$
so from \eqref{stabU} we yield
\begin{equation}
\label{Wv}
\int_{\{\frac{1}{2}<x_n, |s|<\frac{1}{2}\}}
|\partial^{\alpha} W(;\lambda,t_2)|^2\leq
C(\varepsilon^{2\kappa_1}+\varepsilon^{2\kappa_1\mu}
\|v\|^2_{L^2((-T,T);H^1(\Omega))})
\end{equation}
where we assumed that $|t_2|<T_1<T-1$, let $\frac{1}{2}=2\frac{ln (1+2\delta_1)}{\gamma}$ (or $\delta_1=
\frac{e^{\frac{\gamma}{4}}-1}{2}$) and used that because of this choice
$\{\frac{1}{2}<x_n, |s|<\frac{1}{2}\} \subset Q^s(2\delta_1)$.
After choosing $\delta_1$ we will comment on the choice of $\beta, \mu, \kappa_1, T_1, T$ to satisfy \eqref{betakappa}. First we choose $\beta, \kappa_1$ to satisfy the first inequality
\eqref{betakappa} (for example,
$\beta=\frac{\delta_1}{3},
\kappa_1= min (\frac{1}{2}, \frac{\delta_1}{3})$. Then we select $T_1=T_0+\frac{1}{2}$ (to fit the later condition at the end of the proof in the case a)). 
After that we choose (large $T$) to satisfy the third inequality \eqref{betakappa} and then
(small positive) $\mu$ to satisfy the second inequality.
To satisfy the condition 
$C <\tau$ in \eqref{Carleman}
with our choice of $\tau$ in
\eqref{taulambda} we can assume that $\varepsilon <
\frac{1}{C}$, otherwise the bounds \eqref{stability2space},
\eqref{stability2} are obvious.

Integrating \eqref{Wv} over
$\{t_2: |t_2-t_1|<\frac{1}{2}\}$, and using \eqref{stabU0} we obtain
\begin{equation}
\label{U00}
\int_{\{\frac{1}{2}<x_n\}}
|\partial^{\alpha} W(x,0;\lambda,t_1)|^2 dx \leq
C(\varepsilon^{2\kappa_1}+
\varepsilon^{2\kappa_1\mu}
\|v\|^2_{L^2((-T,T);H^1(\Omega))}),
\end{equation}
provided $|t_1|<T_1-\frac{1}{2}$.

Let $\theta=\frac{1}{\lambda}$.
As known, e.g. \cite{R2}, Lemma 6, for solutions $H(t_1,\theta)$ to the heat equation
$\partial_{\theta}H-
\partial^2_{t_1}H=0$ we have
\begin{equation}
\label{H}
\|H( ,0)-H( ,\theta)\|_{(0)}
(\mathbb{R})
 \leq
C \sqrt{\theta}\|H( ,0)\|_{(1)}
(\mathbb{R}).
\end{equation}
Let $\chi_3$ be a $C^{\infty}(
\mathbb{R})$-function,
$\chi_3(t)=1$ if $|t|<T_0$,
$\chi_3(t)=0$ if $T_0+\frac{1}{2}<|t|$,
$0\leq\chi_3\leq 1$. Due to our choice of  $\delta_1$ (depending on $\Omega, \Omega_0, A, a_0,\Gamma$),  $\Omega_0\subset \{\frac{1}{2}<x_n\}$. Using that
$H(t_1,\theta;x)=\partial^{\alpha}W(x,0;\theta^{-1},t_1)$ solves the heat equation and 
$H(t_1,0;x)=\partial^{\alpha}(\chi_{T}(t_1)w(x,t_1))$ we yield
$$
\|\partial^{\alpha}v\|_{(0)}((\Omega_0\setminus \bar\omega)\times(-T_0,T_0))=\|\partial^{\alpha}w\|_{(0)}((\Omega_0\setminus\bar\omega)\times(-T_0,T_0))\leq
\|\chi_3\partial^{\alpha}(\chi_{T}w)\|_{(0)}((\Omega_0\setminus\bar\omega)\times \mathbb{R})\leq
$$
$$
\|\chi_3(\partial^{\alpha}W( ,0;\lambda,0)-
\partial^{\alpha}(\chi_{T}w))\|_{(0)}((\Omega_0\setminus\bar\omega)\times \mathbb{R})+
\|\chi_3(\partial^{\alpha}W( ,0;\lambda, )\|_{(0)}((\Omega_0\setminus\bar\omega)\times \mathbb{R})\leq
$$
$$
\|\partial^{\alpha}W( ,0;\lambda,0)-
\partial^{\alpha}(\chi_{T}w)\|_{(0)}((\Omega_0\setminus\bar\omega)\times \mathbb{R})+
\|\chi_3(\partial^{\alpha}W( ,0;\lambda, )\|_{(0)}((\Omega_0\setminus\bar\omega)\times \mathbb{R})
$$
provided $T_0= T_1-\frac{1}{2}, |\alpha|=1, \alpha_0=0$.
If $\alpha_0=1$ the same argument holds when we substitute $\partial^{\alpha}
(\chi_{T}v)$ by
$i\partial_t
(\chi_{T}v)$.
Using \eqref{U00} and \eqref{H} we conclude that
$$
\|\partial^{\alpha}v\|_{(0)}((\Omega_0\setminus \omega)\times(-T_0,T_0))\leq
$$
$$
\frac{C}{\sqrt{\lambda}}
(\|\partial^{\alpha}(\chi_{T}w)\|_{(0)}(\Omega\times
\mathbb{R})+\|\partial_t\partial^{\alpha}(\chi_{T}w)\|_{(0)}(\Omega\times
\mathbb{R}))+ C(\varepsilon^{\kappa_1}+
\varepsilon^{\kappa_1\mu}
\|v\|_{L^2((-T,T);H^1(\Omega))})
$$
provided $T_0= T_1-\frac{1}{2}\leq T-1-\frac{1}{2}$.
Recalling \eqref{w},\eqref{taulambda} we obtain \eqref{stability2}
and complete the proof of Theorem 1.3 in the case a).

The proof in the case b) is similar. We only indicate  few needed changes.
Instead of \eqref{PDEw} we will have the partial differential equation and the Cauchy data
$$
(A-a_0\partial_t -\partial_t^2) w= f_0,\;f_0=\chi_0 f+
2\sum_{j,m=1}^n a_{jm}\partial_j\chi_0\partial_m v+(A\chi_0-a \chi_0)v\;\mbox{in}\; Q=\Omega\times(-T,T).
$$
\begin{equation}
{\label{CauchyHw}}
   w= \chi_0 v_0,  \partial_{\nu} w=\chi_0 v_1
   +\partial_{\nu}\chi_0 v_0\;  \text{ on }\; \Gamma_0\times(-T,T).
\end{equation}
and accordingly instead of \eqref{AW} 
$$
AW+a_0i \partial_s W+\partial_s^2 W= F_0 -F( ;T)\;\textit{on}\;Q^s,
$$
\begin{equation}
{\label{CauchyE}}
   W= \chi_0 V_0,  \partial_{\nu} W=\chi_0 V_1+\partial_{\nu}\chi_0 V_0\;  \text{ on }\; \Gamma_0\times(-1,1).
\end{equation}
So an application of the Carleman estimate \eqref{Carleman} will give
$$
\int_{Q^s(2\delta_1)} \tau^{3-2|\alpha|}|\partial^{\alpha}W|^2 e^{2\tau\varphi} \leq
C(\int_{Q^s}|F|^2 e^{2\tau\varphi}+
\int_{\Gamma\times(-1,1)}(\tau^3 |V_0|^2+\tau(|\nabla_{x,s} V_0|^2+|V_1|^2) e^{2\tau\varphi})+
$$
$$
\int_{Q^s}|F( ;T)|^2 e^{2\tau\varphi} +\int_{ Q^s\setminus Q^s(\delta_1)}(|W|^2+|\nabla_{x,s} W|^2) e^{2\tau\varphi}, |\alpha|\leq 1.
$$
and consequently
$$
e^{4\tau\delta_1}
\int_{Q^s(2\delta_1)} \tau^{3-2|\alpha|}|\partial^{\alpha}W|^2  \leq
C(e^{2\tau(\Phi+1)}(\|F^2\|_{(0)}(Q^s) +
$$
$$
\tau^3 \| V_0\|^2_{(1)}(\Gamma\times(-1,1))+\tau
\|V_1\|_{(0)}^2(\Gamma\times(-1,1)))+e^{2\tau\Phi}\|F( ;T)\|_{(0)}^2 (Q^s)+
e^{2\tau \delta_1}\|W\|_{(1)}^2(Q^s)).
$$
By using \eqref{Uu} with $\Gamma\times(-T,T)$ instead of $Q$
as in the case a) we arrive at \eqref{stabU}, and the rest of the proof in the case b) is the same as in the case a).

The proof is complete.

\section{Conclusion}

We demonstrated that the solution of the Cauchy or continuation problem is improving with growing wave number  disregard of geometry of a domain
or (non)trapping properties of the metric corresponding to the principal part of the elliptic equation. These results suggest much better
controllability of the higher frequency waves and resolution in the inverse problems (for sources, obstacles, or medium properties) for
any site of (boundary) controls or sensors. 
  
We think that the results of this paper can be extended onto higher order elliptic equations and systems. So far a conditional H\"older stability
is obtained for hyperbolic principally diagonal systems of second order (including isotropic elasticity and Maxwell systems) under (pseudo)convexity conditions 
and sharp uniqueness of the continuation results are obtained without these conditions \cite{EINT}. 
Despite recent results \cite{BKL}, \cite{V} stability bounds in the whole sharp uniqueness domain (given by Tataru \cite{T}) are still not available even in the scalar case. 
An important question is about minimal a priori constraints on the high frequency part of a solution. It is feasible that semi norms $|||\cdot|||_{(m;k)}(\Omega), m=2$ 
in Theorem 1.2 can be replaced by a similar semi norm with $m=1$, imposing only natural energy constraints on the high frequency part of $u$.  Another useful confirmation 
of increasing stability can be obtained by proving that there are growing invariant subspaces where the
solution of the Cauchy problem \eqref{PD}, \eqref{CauchyD} is Lipschitz stable. We will give one of corresponding conjectures.

 Let $\Omega$ be a Lipschitz bounded domain. Let us assume that there is an unique solution $u\in H^{(1)}(\Omega)$ of the following Neumann problem
 $$
 Au+cku+k^2u=0\;\mbox{in}\;\Omega,
 $$
$$
\partial_{\nu}u=0\;\mbox{on}\;\partial\Omega\setminus \Gamma_1,\;
\partial_{\nu}u=g\in H^{(-\frac{1}{2})}(\Gamma_1),\;\mbox{on}\;\Gamma_1.
$$
The operator $B$ mapping $g$ into $u_0=u$ on $\Gamma_0$ is compact from $L^2(\Gamma_1)$ into $L^2(\Gamma_0)$. Hence it admits the singular value decomposition consisting of  
a complete orthonormal system of functions $g_m, m=1,2,...$ in $L^2(\Gamma_1)$ and corresponding singular values $\sigma_m\geq \sigma_{m+1}>0$ (eigenfunction and square roots 
of eigenvalues of $B^*B$). The conjecture is that there are positive numbers $\delta_1,\delta_2$ depending only on $A,c$ and $\Omega$ (but not on $k$) such that $\sigma_m>\delta_1$ 
when $m<\delta_2 k$. This conjecture for some interesting plane $\Omega$ was numerically confirmed in \cite{IK}.

 Use of a low frequency zone does not need any (pseudo)convexity type assumptions on $\Omega, \Gamma, A$ and for this reason is very promising for applications since the observation 
sites $\Gamma$ or $\omega$ can be chosen arbitrarily at convenient locations. In the  paper \cite{IK} we studied this phenomenon on more detail and gave regularization schemes for 
numerical solution incorporating the increasing stability. We gave several numerical examples of increasing
 stability for the Helmholtz equation in some interesting plane domains, admitting or not admitting explicit analytical solution and complete analytic justification. It is important 
to collect numerical evidence of the increasing stability for more complicated geometries and for systems. 

An increasing stability is expected in the inverse source problem, where one looks for $f$ in the Helmholtz equation $(\Delta+k^2)u=f$ 
(not depending on $k$) in $\Omega=\{x: 1<|x|<R\}$ from the Cauchy data $u, \partial_{\nu}u$ on $\Gamma=\{x:|x|=1\}$, $k_*<k<k^*$. General uniqueness results and
convincing numerical examples of increasing stability are given in \cite{IL1}. One needs to obtain stability estimates improving with growing $k^*$ and to give more of 
numerical evidence of better resolution for larger $k^*$. Such stability estimates and numerical results under (pseudo)convexity conditions were obtained in \cite{CIL},
\cite{IL} by using sharp bounds on analytic continuation to higher wave numbers and exact observability for corresponding hyperbolic equations.

It was (numerically) observed, that use of only low frequency zone can produce a stable solution of the inverse problem, where one looks for a speed of the propagation from 
all possible boundary measurements. As above this low frequency zone grows with the wave number. Currently, there are no  analytic proofs of it. One can look at the 
linearized problem: find $f$ (supported in $\Omega\subset \mathbb R^3$) from
$$
\int_{\Omega} f(y) \frac{e^{ki|x-y|}}{|x-y|}\frac{e^{ki|z-y|}}{|z-y|}dy
$$
given for $x,z\in \Gamma \subset \partial \Omega$.  The closest analytic results on improving stability are obtained in \cite{I4}, \cite{IW} for 
the Schr\"odinger potential and in \cite{ILW} for the attenuation and conductivity coefficients. The analytic results
in \cite{ILW} suggest a possibility of better resolution in the electrical impedance tomography when using the data from electromagnetic waves of higher frequencies.

{\bf Aknowledgement}

This research is supported in part by the Emylou Keith and Betty Dutcher Distinguished Professorship and the NSF grant DMS 15-14886.

\end{document}